\documentclass[12pt,leqno]{amsart}

\usepackage{amssymb}
\usepackage{latexsym}

\newcommand\dspace{\lineskip=2pt\baselineskip=18pt\lineskiplimit=0pt}
\dspace

\newcounter{deficislo}[section]
\newcounter{propcislo}[section]
\newcounter{lemacislo}[section]
\newcounter{theocislo}[section]
\newcommand{\defi}[1]{\refstepcounter{deficislo}{\noindent \bf
Definition \thesection.\thedeficislo.\  }{\it #1}}
\newcommand{\pro}[1]{\refstepcounter{propcislo}{\noindent \bf
Proposition \thesection.\thepropcislo.\  }{\it #1}}

\newcommand{\theo}[1]{\refstepcounter{theocislo}{\noindent \bf
Theorem \thesection.\thetheocislo.\  }{\it #1}}

\newcommand\gjk{\gamma_{jk}}
\newcommand\glk{\gamma_{lk}}

\newcommand{\de}{\delta}
\newcommand{\ab}[1]{\vert z\vert^{#1}}
\newcommand{\cdva}{{\mathbb C^2}}

\newcommand{\ml}{M^{k,l}_a}
\newcommand{\te}{\theta}

\newcommand{\zz}{(z,\bar z)}
\newcommand{ \nl}{\newline}

\newcommand{\fzs}{{\textstyle{F^*(z^*,\bar z^*,u^*)}}}

\newcommand{ \al}{\alpha}

\newcommand{\la}{\lambda}
\newcommand{\La}{\Lambda}

\begin{document}
\title{
 Higher order invariants of Levi degenerate
hypersurfaces}
\author {Martin Kol\'a\v r}
\address {Department of Mathematical Analysis, Masaryk University,
Janackovo nam. 2a, 662 95 Brno } \email {mkolar@math.muni.cz }
\thanks{The author was supported by a grant of the GA \v CR no.
201/05/2117 }

\dedicatory{\it Dedicated to Professor Joseph J. Kohn
on the occasion of his 75-th birthday.}

\begin{abstract}
The first part of this paper considers higher order CR invariants
of three dimensional hypersurfaces of finite type. Using a full
normal form we give a complete characterization of hypersurfaces
with trivial local automorphism group, and analogous results
 for
finite groups. The second part
 considers
 hypersurfaces of finite   Catlin
 multitype,
and the Kohn-Nirenberg phenomenon in higher dimensions. We give a
necessary condition for local convexifiability of  a class of
pseudoconvex hypersurfaces in $\mathbb C^{n+1}$.
\end{abstract}
\maketitle

\section{Introduction}

In complex dimension two, the lowest order local CR invariant of a
smooth hypersurfaces  $M \subseteq \mathbb C^2$ at a point $p\in
M$ is the type of the point. This fundamental invariant was
introduced by J. J. Kohn in \cite{K}, in order to study  boundary
behavior of the $\bar
\partial$ operator on weakly pseudoconvex domains.  It  measures
 the number of commutators of CR and anti-CR vector
fields needed to span the whole complexified tangent space to $M$
at $p$. In more geometric terms, the type of $p$ is the maximal
order of contact between $M$ and complex curves passing through
$p$.

On the next level, there is a well defined model hypersurface at
$p$, which gives  important numerical local invariants,
characterizing local geometry of $M$  (\cite{KN}).

In higher dimensions, local  geometry of Levi degenerate
hypersurfaces
 is
much more complicated. The order of contact with complex curves
does not give an open condition anymore, and  cannot characterize
subellipticity of the $\bar{\partial}$ operator on the boundary.
 In order to obtain
invariants relevant for analysis of the $\bar{\partial}$
equation,
one has to consider orders of contact with singular complex
varieties. Such invariants have been introduced in the work of J.
P. D'Angelo (\cite{D}).

Let $M$ be a smooth hypersurface in $\mathbb C^{n+1}$, where
$n\geq 2$. If $d_k$ denotes the maximal order of contact of $M$
with complex varieties of dimension $k$ at the given point, the
n-tuple $(d_n,\dots, d_1)$ is called the D'Angelo multitype of
$p$.

For pseudoconvex hypersurfaces, D. Catlin (\cite{C}) introduced
a different, more algebraic notion of multitype. As an  important
advantage, it provides a well defined weighted-homogeneous model,
an essential tool for local analysis. The two notions of multitype
coincide on a class of hypersurfaces called semiregular
(\cite{DH}), or h-extendible (\cite{Y}). It contains, for example,
all decoupled and  all convexifiable hypersurfaces.

Since the definition of multitype  is nonconstructive, and the
models are not uniquely defined, it is not a priori clear what is
the relation between various models. In some  situations, when low
order boundary invariants are needed,
 it is enough to choose an arbitrary model.
On the other hand, in order to study higher order CR invariants it
is essential to understand the non-uniqueness in the definition of
models. In particular, it is not a priori obvious whether all
models are necessarily  biholomorphically equivalent. In the case
of  h-extendible pseudoconvex hypersurfaces, biholomorphic
equivalence of models was proved by N. Nikolov in \cite{Nik}.

The first part of this paper considers local CR invariants of
hypersurfaces of finite type in $\mathbb C^2$. In dimension two, a
complete set of local CR invariants  can be constructed, on the
level of formal power series, using analysis of generalized
Chern-Moser operators (\cite{Ko1}).
 Combined with the result of M.
S. Baouendi, P. Ebenfelt and L. P. Rothschild on convergence of
formal equivalences, it gives a solution to the local equivalence
problem. One of the main applications of this result provides
understanding of local symmetries of finite type hypersurfaces.

 In Section 3 we
give a characterization of finite type hypersurfaces with trivial
local automorphism group, in terms of normal form coordinates.
Analogous complete classification is obtain in the case of finite
local automorphism group of order $m >1$.

In Section 4 we consider a constructive approach to the Catlin
multitype on a general (not necessarily pseudoconvex) hypersurface
in $\mathbb C^{n+1}$. Section 5 gives a review of Kohn-Nirenberg
invariants in complex dimension two. The results are used in
Section 6, where  we define certain numerical invariants on
hypersurfaces for which at the given point all multitype entries
coincide.
We show that they can
be used to study the Kohn-Nirenberg phenomenon in higher
dimensions. In particular, we obtain a necessary condition for
local convexifiability of such hypersurfaces.

\section{Preliminaries}
In the first part of this paper, we will consider  a real analytic
hypersurface $M\subseteq \mathbb C^2$ and a point $p\in M$ of
finite type $k$.

 In order to  describe $M$ in a neighbourhood of  $p$, we will
 use local
holomorphic coordinates $(z,w)$ centered at $p$, where  $z = x +
iy,\  w = u + iv$. The hyperplane $\{ v=0 \}$ will be assumed
tangent to $M$ at $p$. In this setting, $M$ is described near $p$
as the graph of a uniquely determined real valued function
$$ v = \Psi(z, \bar z, u).$$
 Recall that $p \in
M$ is a point of finite type if and only if there exist local
holomorphic coordinates such that $M$ is given by
\begin{equation}v = P(z,\bar z) + o(\ab k + \vert u \vert ),\ \label{2.9}\end{equation} where the leading
term  is a  nonzero real valued homogeneous polynomial of
degree $k$ without harmonic terms,
\begin{equation}P(z, \bar z) = \sum_{j=e}^{k-e}
a_j z^j\bar z^{k-j} ,\ \label{2.101}\end{equation} with $1 \leq e
\leq \frac{k}2$ and
 \begin{equation}a_e = 1.
\label{ae}
\end{equation}
 We will now recall some basic facts from the normal form construction
\cite{Ko1}.

Using (\ref{2.101}), we will define  two basic
integer valued invariants. The first one,
denoted by $e$ and defined implicitly by (\ref{2.101}), is the
essential type of the model hypersurface (defined below).

When $e < \frac{k}2$, the second invariant is defined as follows.
Let $e=m_0 <m_1 <\dots <m_s <\frac k2$ be
 the indices in (\ref{2.101}) for which $a_{m_i}\neq 0$. The invariant,
 denoted
 by $d$,
 is the greatest common divisor of the numbers
$\ k-2m_0, k-2m_1,  \dots, k-2m_s$.

 The polynomial $P$  need not be determined
uniquely by (\ref{2.101}), (\ref{ae}).
In order to make it unique, the following condition is imposed.
 Denote
$$q_i = \frac {\gcd( k-2m_0, k-2m_1,  \dots, k-2m_i)}{
\gcd ( k-2m_0, k-2m_1,  \dots, k-2m_{i+1})}$$
 for $0 \leq i \leq
s-1$. In addition to (\ref{2.101}), (\ref{ae}),  we require that $P$
satisfies
\begin{equation}
 \arg a_{m_{i+1}} \in [0,\frac{2\pi}{q_i})\ \label{arg}
\end{equation}
for $0 \leq i \leq s-1$.
 This  determines  $P$ uniquely.

 The model hypersurface $M_H$
to $M$ at $p$ is defined using the normalized leading homogeneous
term,

\begin{equation} M_H = \{(z,w) \in \cdva\ | \
 v  = \sum_{j=e}^{k-e}
a_j z^j\bar z^{k-j} \}. \label{2.11}\end{equation}
 In particular, when the leading term is circular, we  write
\begin{equation}S_k = \{(z,w) \in \mathbb C^2 \ \vert \   v = \ab k
\}. \end{equation}
 Another exceptional model is the tubular hypersurface
\begin{equation}T_k = \{(z,w) \in \Bbb C^2 \ \vert \
v = \frac1k [(z+\bar z)^k - 2Re\;z^k] \}. \end{equation}

 The local automorphism group of a hypersurface
 $M$ at a point $p$
will be denoted by $Aut(M,p)$. It was proved in \cite{Ko1} that if
$e < \frac{k}2$, the local automorphism group of $M_H$ consists of
transformations
$$z^* =   \delta \exp {i\theta} z,\ \  \ \ \   w^* = \delta^k w,$$
                           where
$\exp {i\te}$ is a $d$-th root  of unity and $\delta > 0$ for $k$
even or $\delta \in \Bbb R \setminus \{ 0 \}$ for $k$ odd.
Hence $Aut(M_H, p) = \Bbb R^+
\oplus {\Bbb Z}_d$ for $k$ even and $Aut(M_H, p) = \Bbb R^* \oplus
{\Bbb Z}_d$ for $k$ odd.

The local automorphism group of $S_k$ is three dimensional,
consisting of transformations of the form

\begin{equation}f( z,w)  =  \frac{ \delta
e^{i\theta} z}{(1 + \mu w)^{\frac1e}}, \ \ \ \ g( z,w)=
  \frac{ \delta^k w}{1 + \mu w},\label{4.1}\end{equation}
                           with $\delta > 0,$ and $\te, \mu  \in \Bbb R$.

\section{Local automorphism groups in $\mathbb C^2$ }
We write
$$\Psi(z,\bar z, u) = P(z,\bar z ) + F(z, \bar z, u),$$
where
$$ F(z, \bar z,u) = \sum_{j,l,m} a_{jlm} z^j \bar z^l u^m. $$
 We will also  consider the Taylor expansion of $F$ in terms of $z,
\bar z$,
$$F(z, \bar z, u) = \sum_{j,l}F_{jl}(u)z^j\bar z^l,$$
where
 $$F_{jl}(u) = \sum_{m} a_{jlm} u^m.$$
The results of \cite{Ko1} give three different complete normal
forms, depending on the form of the model. There are  two
exceptional models, $S_k$ and $T_k$, while the generic case covers
all remaining models.

 When $e = \frac{k}2$, we have
the following complete  normalization.
  $F$ is  in  normal form if
  \begin{equation}
\begin{array}{rl} F_{j0} &= 0, \ \ \ \ \ j=0,1,\dots,  \\
F_{e,e+j} &= 0, \ \ \ \ \ j= 0,1,\dots, \\
F_{2e, 2e} & = 0,\\
F_{3e, 3e} & = 0,   \\
F_{2e, 2e-1}& = 0. \label{fo1} \end{array} \end{equation}

Normal coordinates
 (i.e.
such in which the normal form conditions hold), are determined
uniquely up to the action of the symmetry group (\ref{4.1}).

When $M_H = T_k$, we have the following normal form conditions:
\begin{equation}
\begin{aligned} F_{j0} &= 0, \ \ \ \ \ j=1,2,\dots,  \\
F_{k-1+j,1} &= 0, \ \ \ \ \ j= 0,1,\dots, \end{aligned}
\label{f2a} \end{equation} and
\begin{equation} F_{2k-2, 2} = Re \; F_{k-2,1} = Re \; F_{k, k-1}  = 0.
 \label{f2}
 \end{equation}

Again, normal coordinates are determined uniquely up to the action
of the symmetry group $Aut(T_k, 0)$.

 Now, let $M_H$ be a generic
model, i.e. $e < \frac{k}2$ and $M_H$ is different from $T_k$.
Denote $F_{k-1}(u) = (F_{1, k-2}(u),
 F_{2, k-3}(u), \dots,
F_{k-2,1}(u))$.
 The normal form conditions
 are:
 \begin{equation}\begin{array} {rl}
F_{j0} & = 0, \ \ \ \ \ j=1,2,\dots,  \\
F_{k-e+j,e} & = 0, \ \ \ \ \ j= 0,1,\dots, \\
F_{2k-2e, 2e} & = 0, \\
(F_{k-1}, P_z) & = 0,  \label{f3} \end{array}\end{equation}
 where
\begin{equation}
(F_{k-1}, P_z) = \sum_{j=1}^{k-2}F_{j,k-1-j}  (j+1)\bar
a_{j+1}.
\end{equation}
The corresponding normal coordinates are unique  up to the action
of the symmetry group $Aut(M_H, 0)$.
\\[2mm]
The following result was obtained in  \cite{Ko2}.
\\[2mm]
 \pro{\it Let $M$ be a Levi degenerate hypersurface
 of finite type
with $e=\frac{k}2$,  not equivalent to $v = |z|^k$. Then
 all local automorphisms expressed in normal
coordinates have the form of  decoupled linear transformations
\begin{equation}
z^* = \de e^{i \theta} z, \ \ \ \ \ w^* = \delta^k w,
\label{dec}\end{equation} for some  $\theta \in \mathbb R$ and
$\delta \in \mathbb R^*$.}
\\[2mm]

Now we turn to the general case. Let $(z,w)$ be normal coordinates
for $M$ at $p$, and $v = P(z,\bar z)  + F(z, \bar z, u)$ be the
defining equation in such coordinates. Denote
$$\Theta_1 = \{(j,l,m) \in \mathbb Z_+^3 :\  j\neq l \ \
\text{and} \ \ \ a_{j,l,m} \neq 0 \}$$ and
$$\Theta_2 = \{(j,k-j,0) \in \mathbb Z_+^3 :\  j\neq k-j \ \
\text{and} \ \ \ a_{j} \neq 0 \}.$$ Set $\Theta = \Theta_1 \cup
\Theta_2$.
 If $\Theta $ is
nonempty, we define
$$\mu_0 = {\gcd}_{(j,l,m) \in \Theta} \ |j-l|.$$

We have the following complete description of all finite type
hypersurfaces with finite stability group.
\\[2mm]
 \theo{ Let $M$ be a Levi degenerate hypersurface
 of finite type, not equivalent to a model hypersurface.
 Then  $Aut (M,p)$ is finite if and only if $\Theta$ is nonempty.
 In this case, the stability group is isomorphic  to $\mathbb
 Z_{\mu_0}$.}
\\[2mm]
 {\it proof:}
First we will prove that in normal coordinates all local
automorphisms of $M$ are decoupled linear, of the
 form (\ref{dec}).
By Proposition 3.1, it remains to consider the case $e<\frac k2$.

By definition, any local automorphism of $M$ in normal coordinates preserves
normal form. In \cite{Ko1}, such transformations are completely
characterized, and correspond to the action of the  local symmetry
group of the model.

Starting with the generic case, we assume $M_H$ is a generic
model, and consider the normal form conditions (\ref{f3}).
 The local symmetry
group of the model acts on normal forms. We will prove that the
transformation by  each element of $Aut(M_H,p)$ preserves the
normal form, hence its action on normal forms is direct and no
renormalization is needed. Since every element of $Aut(M_H,p)$ is
a decoupled linear transformation, its application clearly
preserves the first three conditions. In order to see that the
last condition is also preserved, we write $P$ as
$$ P( z, \bar z ) = \sum_{j=1}^{k-1} a_{j} z^{j}\bar z^{k-j},$$
where $ a_j \neq 0$ only if $k-2j$ is divisible by $d$.
Denote
$$ \beta_j(u) = F_{j,k-1-j}(u).$$
We have
$$ P_z(z, \bar z) =\sum_{j=1}^{k-1} j a_{j} z^{j-1}\bar
z^{k-j}. $$
 Since dilations clearly preserve all the normal form conditions, we
  consider the action of a transformation $z^* = \al z$,
where $\al^d = 1.$ Using this and the fact that $a_j \neq 0$ implies
$k-2j$ is divisible by $d$, we obtain

$$(F^*_{k-1}, P_z) = \sum_{j=2}^{k-1}
j \bar a_{j}  \beta_{j-1}(u) \al^{j-1} \bar \al^{k-j}
  =
  $$
  $$=\al^{-1} \sum_{j=2}^{k-1}
j \bar a_{j}  \beta_{j-1}(u) \al^j  \bar \al^{k-j}  = \al^{-1}(F_{k-1},
P_z) = 0,$$
which proves the claim.
Now let the model be the tubular hypersurface $T_k$.
In this case $e=1$,  and since all the coefficients $a_j$ of $P$
are nonzero, we obtain immediately the value of $d$. If $k$ is
even, all the numbers $\ k-2m_0, k-2m_1,  \dots, k-2m_s$ are even,
and we have $d=2$. If $k$ is odd, then again immediately $d=1$.

Consider  the normal form conditions (\ref{f2a}), (\ref{f2}).
Clearly this normalization is preserved when an element of
$Aut(T_k,0)$ is applied.

Thus we have proved that whenever $M$ is different from $S_k$, all
local automorphisms  in normal coordinates are of the form
(\ref{dec}). Hence it remains to consider the action of
(\ref{dec}) on the defining equation of $M$.

If $\Theta $ is empty, then any rotation in $z$ preserves $M$, so
the local automorphism group is infinite. Let us assume now that
$\Theta$ is not empty.

Since all local automorphisms are decoupled linear in normal
coordinates, they act on each monomial in the expansion of $F$
separately, in an obvious way. Since $M$ is not a model, it
follows immediately that weighted dilations do not preserve $F$.
If $(j,l,m) \in \Theta$, then the coefficient $a_{j,l,m}$ is
preserved by a rotation $z^* = \exp i\theta \; z$ if and only if
$$(\exp i\theta)^{j-l} =1.$$
 Since this holds for any element of $\Theta$,  the result follows.
\\[2mm]
As a particular case, we obtain a complete description of
hypersurfaces with trivial local automorphism group.
\\[2mm]
 \theo{ Let $M$ be a Levi degenerate hypersurface
 of finite type, not equivalent to a model.  $Aut (M,p)$ is trivial if and only if $\Theta$
 is nonempty and $\mu_0 = 1$.}

\section{Hypersurfaces of finite multitype}
Let  $M \subseteq \Bbb C^{n+1}$ be a smooth hypersurface (not
necessarily pseudoconvex),  and $p $ be a Levi degenerate point on
$M$. We will assume that $p$ is a point of finite type in the sense
of Bloom and Graham.
In this section we consider a constructive approach to the
Catlin's definition of multitype.

 Consider  local holomorphic coordinates
$(z,w)$, where $z =(z_1, z_2, ..., z_n)$ and $w=u+iv$, $z_j = x_j
+ iy_j$, centered at the point $p \in M$. Again, the hyperplane
$\{ v=0 \}$ is assumed to be tangent to $M$ at $p$. $M$ is
described near $p$ as the graph of a uniquely determined real
valued function
\begin{equation} v = \Psi(z_1,\dots, z_n,  \bar z_1,\dots,\bar z_n,  u).
\label{vp}
\end{equation}
We now apply Catlin's definition of multitype to $M$ at $p$. In
the following, $\al, \beta $ will denote multiindices, and we
will use the standard multiindex notation.
\\[2mm]
\defi{A weight is an n-tuple of nonnegative
 rational numbers $\La = (\la_1, ...,
\la_n)$, where $0 \leq\la_j\leq \frac12$, and $\la_j \ge
\la_{j+1}$, such that for each $k$ there exist nonnegative integers $l_1, ...,
l_n$ satisfying $l_k > 0$ and
$$ \sum_{j=1}^n l_j \la_j = 1.$$
}
\\[2mm]
The component $\la_j$ of $\Lambda$ is interpreted as the weight of the
variable $z_j$.
The variables $w$ and $u$ are given weight one.
 The weighted degree of a monomial $c_{\al
\beta l}z^{\al}\bar z^\beta u^{l} $ is
$$ wt(c_{\al \beta l}z^{\al}\bar z^\beta u^{l}) = l +  \sum_{i=1}^n (\al_i + \beta_i ) \la_i.$$
A real valued polynomial $P(z, \bar z, u)$  is $\La$ - homogeneous of
weighted degree $\kappa$ if it is a sum of
 monomials of weight $\kappa$.
\\[2mm]

 A weight $\La$ will be called distinguished if there exist
coordinates $(z,w)$ in which the defining equation has form

\begin{equation} v = P\zz + o_{wt}(1),\label{1}\end{equation} where $P\zz$ is
a $\La$ - homogeneous polynomial of weighted degree one which is
not pluriharmonic, and $o_{wt}(1)$ denotes terms in the Taylor
expansion of weight greater than one.

The fact that distinguished weights do exist follows from the
assumption of Bloom-Graham finite type (\cite{BG}).

We denote by  $\Lambda_M = (\mu_1, \dots, \mu_n)$  the
infimum of  distinguished weights with respect to the lexicographic
ordering.

The multitype of $M$ at $p$ is defined to be the
n-tuple $(m_1, m_2, \dots, m_n)$, where $m_j = \frac1{\mu_j}$ if
$\mu_j \neq 0$ and $m_j = \infty $ if $\mu_j = 0$. If none of the
$m_j$ is infinity, we say that $M$ is of finite multitype at $p$.

Note that since the definition of multitype  considers all distinguished weights, the
infimum is a biholomorphic invariant, and we may speak of {\it
the} multitype.

Coordinates corresponding to a distinguished weight $\Lambda$, in
which the local description of $M$ has form (\ref{1}), with $P$
being  $\Lambda$ - homogeneous,  will be called $\La$ - adapted.

$\Lambda_M $ will be called the multitype weight. Note that for
any $\delta > 0$ there exist only finitely many points $(l_1,
\dots, l_n) \in \mathbb Z^n_+$ such that $l_j \leq \frac1{\delta}$
for all $j = 1, \dots, n$. It follows immediately  that if $M$ at
$p$ is of finite multitype, $\La_M$ - adapted coordinates do
exist (cf. \cite{C}).

From now on we assume that $p \in M$ is of finite multitype.
 If
(\ref{1}) is the defining equation in some $\La_M$ - adapted
coordinates, we define a model hypersurface to $M$ at $p$ to be

\begin{equation} M_H = \{(z,w) \in \mathbb C^{n+1}\ | \
 v  = P \zz \}. \label{2.10}\end{equation}

Models are useful for many geometric and analytic results. In
order to deal with biholomorphisms between models, we introduce
the following terminology. Here weighted degree is understood with
respect to the multitype weight $\La_M$.
\\[2mm]

\defi{
  A transformation
$$ w^* = w + g(z_1, \dots z_n, w),\ \ \ \ \ z_i^* = z_i
+ f_i(z_1, \dots z_n, w)$$ preserving form (\ref{vp}) is called

-- homogeneous if $f_i$ is a $\Lambda_M$-homogeneous polynomial of
weighted degree $\mu_i$ and
 $ g $
is a $\Lambda_M$-homogeneous
polynomial of weighted degree one,

-- subhomogeneous if $f_i$ is a polynomial consisting of monomials
of weighted degree less or equal to $\mu_i$ and $g$ consists of
monomials of weighted degree less or equal to one,

-- superhomogeneous if the Taylor expansion of $f_i$ consists of
terms of weighted degree greater or equal to $\mu_i$ and $g$
 consists of
terms of weighted degree  greater or equal to one.
 }
\\[2mm]

We write $P$ in the form

\begin{equation}
P(z, \bar z) = \sum_{| \al, \al' |_{\La_M} =1} A_{\al, \al'} z^{\al} \bar
z^{\al'},
\end{equation}
where  $|\al, \al'|_{\La_M} = \sum_{j=1}^n \mu_j
(\al_j+{\al'}_j)$.  Homogeneous transformations are of the form

\begin{equation}
z^*_i = z_i + \sum_{|\al|_{\La_M} = \mu_i}C_{\al}z^{\al},
\ \ \ \ \ w^* = c w + \sum_{|\al|_{\La_M} = 1}D_{\al}z^{\al}
\end{equation}
where $|\al|_{\La_M} = \sum_{j=1}^n \mu_j \al_j$ and $c \in
\mathbb R^*$.
\\[2mm]
Let us remark that   the problem of biholomorphic equivalence of
models  is considered in \cite{Ko7}. The following result was
obtained there.  Note that  models are understood in the sense of
this section, i.e. corresponding to $\La_M$-adapted coordinates.
\\[2mm]
\theo{
A biholomorphic transformation
takes $\La_M$-adapted coordinates into $\La_M$-adapted coordinates
if and only if it is superhomogeneous. Moreover, let $M_H$ and
$\tilde M_H$ be two models for $M$ at $p$. Then there is a
homogeneous transformation which maps $M_H$ to $\tilde M_H$. In
particular, any two models are biholomorphic by a polynomial
transformation.}

\section{Kohn-Nirenberg hypersurfaces }

 We now review some
explicit conditions for local convexifiability of pseudoconvex
hypersurfaces in complex dimension two, which will be used in the
next section. There are no new results in this section.

Let
\begin{equation}
P(z, \bar z)= a_0\ab k +   \sum_{j=2,4,...,k}\ab{k-j} Re(a_j z^j)
\label{3} \end{equation} be a subharmonic but not harmonic
homogeneous polynomial of degree $k$. We will denote
$$\glk = \frac{k}{l^2-k} $$
if $l^2\geq 3k-2$ and $$\glk = \sqrt{\frac{(4k - l^2 -
4)k^2}{(4k-4)(k^2 -l^2)}}$$ if $l^2\leq 3k-2$.

Further, we consider hypersurfaces of the Kohn-Nirenberg type.
 Let
$$M^{k,l}_a =
\{(z,w)\in \mathbb C^2\ \mid \ Im \ w = P^{k,l}_a(z, \bar z)\},$$ where
$$P^{k,l}_a \zz = \ab k +a\ab{k-l}Re\ z^l$$
with $a \geq 0$.
\\[2mm]
We have the following characterization of convexifiability of
$M^{k,l}_a$,  obtained in \cite{Ko5}.
\\[2mm]
\pro{ $\ml$ is convex if and only if $a\le \gamma_{lk}$.
 Moreover, if $l$ is
not a divisor of $k$, then this condition is equivalent to
convexifiability of $\ml$. }
\\[2mm]
The following result appears  in \cite{Ko51}, as Theorem 3 .
\\[2mm]
\pro{ Let the model  at $p \in M$ be given by  (\ref{3}). If $M$
is convexifiable at $p$, then \nl (i)$\ \dfrac{ \vert a_j
\vert}{a_0}  \leq\gjk $ for all $j>\frac{k}2,$ \nl and \nl (ii)$\
\dfrac{ \vert a_j \vert}{a_0}  \leq 2 \gjk $ for all $j \leq
\frac{k}2$. }
\\[2mm]

\section{Hypersurfaces with homogeneous models}

 We will now consider a smooth
pseudoconvex hypersurface $M \subseteq
\mathbb C^{n+1}$  and 
local holomorphic coordinates $(z_1, z_2, ..., z_n, w)$, where
$w=u+iv$ and $z_j = x_j + i y_j$, centered at a point $p \in M$. We
assume $p$ is of finite Catlin multitype.

Consider $\La_M$ - adapted coordinates, in which the hypersurface
is described by
$$ v = P(z_1,\dots, z_n,  \bar z_1, \dots, \bar z_n ) + o_{wt}(1).$$

 It is well known that on any
locally convex domain the Catlin and D'Angelo multitypes coincide,

 and the numbers $m_j$ are all even integers
 (see \cite{Yu2}).
 Hence inequality of the two multitypes is a trivial
 obstruction
 to convexifiability. Our aim is to study other possible
 obstructions, hence we restrict ourselves to domains
 on which the two multitypes coincide.
 This simplifies
substantially the form of the leading polynomial $P\zz$. When
restricted to a coordinate axis $z_j$ it gives a subharmonic but
not harmonic real valued homogeneous polynomial of degree $m_j$ of
the form
\begin{equation}P_j(z_j, \bar z_j)= a_0^j\vert z_j \vert^{m_j} +
   \sum_{i=2,4,...,m_j}\vert z_j \vert^{m_j-i} Re(a_i^jz_j^i)
 \label{pj}  \end{equation}
for some $a_i^j \in \mathbb C$ and $a_0^j > 0$.

We will now  consider hypersurfaces  for which the multitype at
the given point satisfies $m_1 = m_2 = \dots = m_n =m$, each entry
being equal to a fixed even integer $m$. Hence, in $\La_M$
- adapted coordinates,  the leading polynomial $P$ is a plurisubharmonic homogeneous
polynomial of degree $m$ which, by the equality of multitypes,
 is
not harmonic along any complex line passing through the origin.
Indeed, if it were harmonic along such a line, the order of
contact with complex curves would exceed $m$.
 In
this case for all homogeneous and subhomogeneous transformations
the $f$ component of the transformation is linear.
This substantially simplifies the analysis.

We define the Kohn-Nirenberg numbers of $M$ at $p$ as follows. For a
nonzero vector $c = (c_1, \dots c_n) \in \mathbb C^n$ we consider
the restriction of $P$ to the complex line
\begin{equation} \Gamma_{c} = \{ z\in \mathbb C^n \ ; \ z = \zeta c, \ \zeta
 \in \mathbb C \},\label{gama}\end{equation}
 generated by $c$. This restriction as a function of $\zeta$ is a subharmonic polynomial
of the form  (\ref{3}), which we denote $P_c$, and its
coefficients by $a_j^{c}$. For an even integer $l$  we define the
Kohn-Nirenberg number
  $$ \kappa^l_M = \sup_{c}\frac{\vert a_l^{c} \vert}{a_0^c}.$$
We have the following necessary condition for local
convexifiability.
\\[2mm]
\pro{ If there exists  an $l>\frac{m}2$ such that
$\kappa^l_M > \gamma_{lm},$
or an $l \leq \frac{m}2$ such that $\kappa^l_M > 2\gamma_{lm},$
then $M$ is not locally convexifiable. }
\\[2mm]
{\it proof:} By Proposition 5.2., $M$ is not convex in the
original coordinates $(z,w)$. Let $(z^*, w^*)$ be another system
of local holomorphic coordinates, and let the biholomorphic
coordinate change be given
 by
\begin{equation}\begin{aligned} z_i^*&=  f_i(z,w) \\
 w^*&=g(z,w). \label{4} \end{aligned} \end{equation}
We may restrict attention to transformations which preserve the
form (\ref{vp}), and moreover satisfy the normalization condition
$g_w(0,0) = 1$. The general case is obtained from this  by an
affine transformation,
 which does not affect convexity.

 Let  $ F^*$
 denote
the function which describes  $M$  in the new coordinates. By
substituting (\ref{4}) into $v^*=\fzs$, we obtain
\begin{equation}\begin{aligned}  F^*(f(z,u+iF(z, \bar z, u)),
\overline{f(z,u+iF(z, \bar z, u))}&,\\ Re\
 g(z,u+iF(z, \bar z, u))   =
  Im\ g(z,Re \ g  &(z,u+iF(z, \bar z,
u)).\label{covf}\end{aligned} \end{equation} Now we will
distinguish two cases. First, let  $g$ contain terms of weight
less than one. Let $g_{\kappa}$ denote the leading homogeneous
term in $g$, where $\kappa$ denotes the corresponding weight. By
comparing  terms of weight $\kappa$ in (\ref{covf}), it follows
that $F^*$ starts with a nonzero pluriharmonic polynomial
$Q(z,\bar z)$ of weight $\kappa$. We choose a complex line
$\Gamma_c$ of the form (\ref{gama}),  such that the restriction of
$Q$ to $\Gamma_c$ is nonzero. The defining equation restricted to
this line has form
$$ v = Re\; \al \zeta^{m \kappa} + o(m \kappa),$$
where $\al \neq 0$,
hence $F^*$ is not locally convex.

Let now  $g$  contain only terms of weight greater or equal to
one, and let $g_1$ denote the homogeneous part of $g$ of weight
one. We separate the leading linear term in $f(z,w)$ in the
$z$-variables. Write
$$f_i(z,w) = L_i(z) + o_{wt} (\frac1m)$$
and denote $L(z) = (L_1(z), \dots, L_{n}(z))$.
 Let  $P^*$  denote the leading homogeneous term in $F^*$. Since
there are no terms of weight less than one on the right hand side
of (\ref{covf}), $P^*$ is of weight one. For terms of weight one
in (\ref{covf}) we obtain
$$P^*(L(z), \overline{L(z)} )= P(z,\bar z) + Im\; g_{1}(z).$$
It follows that  the leading term in $F^*$ is obtained from $P$ by
a linear transformation in $z$ and addition of pluriharmonic
terms. In the original coordinates we choose a line $\Gamma_{c}$,
on which $\frac{\vert a_l^{c} \vert}{a_0^c} > \gamma_{lm},$ for
some $l>\frac{m}2$ or $\frac{\vert a_l^{c} \vert}{a_0^c}
> 2\gamma_{lm},$ for some $l \leq \frac{m}2$, and denote
$\Gamma^*_c$ the image of this line by the linear part $L(z)$. It
follows from linearity that  the same condition holds on this
line, hence by Proposition 5.2., $F^*$ is not convex.

\end{document}